\documentclass[10pt]{amsart}
\usepackage{amsmath,amsthm,amscd,amsfonts,amssymb}
\usepackage[all]{xy}

\begin{document}
\newtheorem{defn}{Definition}[section]
\theoremstyle{definition}
\newtheorem{thm}[defn]{Theorem}
\newtheorem{propo}[defn]{Proposition}
\newtheorem{cor}[defn]{Corollary}
\newtheorem{lem}[defn]{Lemma}
\theoremstyle{remark}
\newtheorem{rem}{Remark}
%\theoremstyle{example}
%\newtheorem{exam}{Example}
%\theoremstyle{proof}
%\newtheorem{pf}{Proof}

%%%%%% 1.- OPERADORES

\newcommand\aut{\operatorname{Aut}}
\newcommand\g {{\operatorname{G_\lambda}}(O_D^n)}
\newcommand\gi {{\operatorname{G_{\lambda+i}}}(O_D^n)}
%%%%%% 2.- CALIGRAFIA
\renewcommand\o{{{\cal O}}}      
\newcommand\Z{{\Bbb Z}}    
\newcommand\s{\sigma}
\newcommand\w{\widetilde}
%%%%%% 3.- ABREVIATURAS
\newcommand\cal{\mathcal}
\newcommand\fu{\underline}
\newcommand\limi{\varinjlim}
\newcommand\limp{\varprojlim}
\newcommand\f{\underset{{\Bbb F}_q}{\otimes}}

%%%%%%%%%%%%%%%%%%%%%%%%%%%%%%%%%%%%%%%%%%%%%%%%%%%%%%%%%%%%%
\title [ Uniformizers for elliptic sheaves ]{  Uniformizers for elliptic sheaves
 }

\author[A. \'Alvarez]{A. \'Alvarez${}^*$}

\address{Departamento de Matem\'atica Pura y Aplicada \\
Universidad de Salamanca \\ Plaza de la Merced 1-4 \\
Salamanca 37008. Spain.}

\thanks{This work is partially supported by the CIGYT research
contract n. PB91-0188 \\$*$  Departamento de Matem\'atica Pura y Aplicada. 
Universidad de Salamanca}

\date{September 1998, MSC-11G09}

\maketitle 

%%%%%%%%%%%%%%%%%%%%%%%%%%%%%%%%%%%%%%%%%%%%%%%%%%%%%%%%%%%%%

\tableofcontents

%%%%%%%%%%%%%%%%%%%%%%%%%%%%%%%%%%%%%%%%%%%%%%%%%%%%%%%%%%%%%

\section{Introduction.}

In this paper we define $k$-elliptic sheaves, $A$-motives and $t$-modules over $A$, which are
obvious generalizations of   elliptic sheaves, $t$-motives and $t$-modules. Following
   results of \cite{An1}, \cite{D}, \cite{LRSt}, \cite{Mu}, \cite{St},... we shall obtain
the equivalence of this objects.
Bearing in mind \cite{Al} we also describe a correspondence between   $k$-elliptic sheaves with
formal level structures (\ref{form}) and discrete subspaces. The moduli scheme for
these objects shall be a subscheme of Sato's infinite Grassmannian, in this setting  the
moduli for
classical $t$-motives will be a closed subscheme of this last subscheme. In the section 6 for
these discrete subspaces, we show a result analogous  to
\cite{An1} for the behaviour of the determinant.

In the same way $A$-motives with formal level structures have associated locally dense subspaces
(\ref{den}) these subspaces are determined by a subspace of formal series of dimension $k$,
called subspace of uniformizers (\ref{unif}). With this we see that  a 
$k$-elliptic sheaf with formal level structure could be seen as a orbit of the action of
$Gl_n({\Bbb F}_q[[t_x]])$ over the set of these subspaces of uniformizers. These subspaces
allow us to get an injective morphism between the moduli scheme of $k$-elliptic sheaves with
formal level structures and  a Grassmannian  of finite $k$-dimensional
subspaces of formal series, Definition \ref{grass}.   An explicit computation can be made
for these subspaces of uniformizers for $k=1$ and one can  check certain analogy between
these subspaces of uniformizers and Dirichlet series for Drinfeld modules
\cite{G}.  These uniformizers allow us  to obtain the arithmetic counterpart of the Baker
function, 
\cite{SW}, defined in the setting of the theory of soliton equations (see \cite{BlSt} for elliptic
sheaves together an   introduction of soliton theory,  an intensive study for $1$-elliptic modules
is done in  \cite{An2}). We can also make explicit the action  of the
classical arithmetics groups, $Gl_n({\Bbb F}_q((t_x)))$, over these subspaces of uniformizers,
with this result we shall obtain the ''elliptic'' Baker function in an explicit way. Furthermore
we can translate  the results of the anti-equivalence of the Krichever functor
\cite{MR}, \cite{Q}, \cite{SW}, to these $1$-dimensional subspaces and  gain similar objects to
Schur pairs for elliptic sheaves, that we will call   elliptic sheaf triples.    As immediate
consequence of these last results, the subspaces of uniformizers also determine the curve where the
elliptic sheaf is defined.  We can
 get   a moduli scheme for these triples. 

As for discrete subspaces, we shall study  the determinant
for locally dense subspaces associated to $A$-motives of $\s$-rank
$1$ and rank $n$ (Drinfeld modules) with formal level structures. We give an open condition in the
moduli scheme of $A$-motives of $\s$-rank
$1$  and rank $n$ with formal level structures, which implies that the determinant of a certain
''scattering'' matrix, \cite{An1}, obtained from the   uniformizers, generates as
$A\otimes
\o_S\{\s\}$-module, the   determinant of the locally dense subspace associated to the $A$-motive.
When
$A={\Bbb F}_q[t]$, we see that this ''scattering'' determinant gives the subspace of uniformizers
associated to the 
$A$-motive determinant. It would be very interesting to obtain a similar result for arbitrary $A$.

 Let us now  briefly state the   contents of the different sections of this paper.
In the second section we define $k$-elliptic sheaves, $A$-motives and $t$-modules over $A$ with
theirs level structures and we recall the ''equivalences'' among these objects, after this
 we are recalling some facts  about Sato's infinite Grassmannian
scheme and vector bundles,  we see 
$k$-elliptic sheaves with formal level structures  could be considered as discrete subspaces
\label{discr} and the moduli functor for this objects is representable by a closed subscheme of 
Sato's infinite Grassmannian scheme.

In the third section we prove
$A$-motives with formal level structures  give in a unique way that subspaces of uniformizers and
we   construct a moduli scheme for these objects.

In the fourth section we consider discrete subspaces attached to elliptic sheaves and the elliptic
Baker function.

In the fifth section we study an  analogous structure to classical Schur pairs in the set of
Drinfeld modules, elliptic sheaf triples.

In the sixth  section  we settle some results about the determinant in our setting.

\bigskip\bigskip

%%%%%%%%%%%%%%%%%%%%%%%%%%%%%%%%%%%%%%%%%%
%%%%%%%%%%%%%%%%%%%%%%

 {\bf List of notation} 

${\Bbb F}_q$      a finite field of $q$ elements ($q=p^m$, $p$  prime)
\bigskip

 $X$         a smooth, proper and geometrically irreducible curve over ${\Bbb F}_q$

\bigskip

$\otimes$   denotes $ \underset {{\Bbb F}_q} \otimes$

\bigskip
$g$          is the genus of $X$

\bigskip
 $\infty$    a  rational point of $X$

\bigskip

$A$      $=H^0(X-\{\infty\},{{\cal O}_X})$

\bigskip

$\Omega_X,\Omega_A$    the sheaves of differentials on $X$ and $A$ respectively

\bigskip

 $x\neq\infty$    a  point of $Spec(A)$

\bigskip

$m_x$ is the maximal ideal of $x$

\bigskip
$t_x$      a uniformizer at $x$

\bigskip

$k(x)$           residue field of $x$
\bigskip

$deg(x)$   =$dim_{{\Bbb F}_q}(k(x))$
\bigskip

$v_x$   is the valuation associated to $x$
\bigskip

$\Omega_x$     is the completion of $\Omega_X$ at $x$

\bigskip
 $\bar{D}$      is an  effective divisor over $Spec(A)$
\bigskip

$I_{\bar{D}}$         is the ideal in $A$ associated to $\bar{D}$
\bigskip

 $S$        an   arbitrary scheme over ${{\Bbb F}_q}$  
\bigskip

$e$       the structure map $S \to {{\Bbb F}_q}$

\bigskip
$\pi$    the natural projection $ X \times S\to S$

\bigskip
$\s$          Frobenius endomorphism over $S$ ($\s (\lambda)=\lambda^q$)

\bigskip
$\s^*(.)$       is the inverse image by $\s$ of the sections on $S$ of a sheaf defined in the
                 category of ${{\Bbb F}_q}$-schemes.

\bigskip

$F$             $=Id\times\s:X\times S \to X\times S$ 
\bigskip

$({\Bbb G}_a)_S$         the additive line group over $S$
\bigskip

%$\w {\Bbb A}^\infty$  ring of adeles outside $\infty$ and 
%$$\w {O}^\infty=\underset {I{\text{ ideal of }}A}   {\underset {I}   \varprojlim}A/I$$

%%%%%%%%%%%%%%%%%%%%%%%%%%%%%%%%%%%%%%%%%%%%%%

\section{$k$-Elliptic sheaves, $A$-motives and $t$-modules over $A$.}

 The results of the beginning of this  section are essentially obtained from
\cite{An1},
\cite{D},
\cite{LRSt},
\cite{Mu}.  Moreover, we can get,  in a easy way, the same results of \cite{Al} for these
objects. 
 Some times to easy notation we are going to consider
$S=Spec(B)$.
$End_B({\Bbb G}_a)_B= B\{\s\}$, where
$B\{\s\}$ is the non-commutative polynomial ring with the commutation rule $(\s.b=b^q\s)$.

\begin{defn}
 An $A$-motive of rank $n$ and $\s$-rank $k$ over $B$ is an $A \otimes B\{\s\}$-module, ${{\cal
N}}$,  locally free of rank $n$ as $A \otimes B$-module, and free of rank $k$ as $B\{\s\}$-module. 
\end{defn}

From this definition we obtain ${\cal N}= {\underset {i\in {\Bbb N}}  \bigcup}{\cal N}_i$,
where
$${\cal N}_i=\{p_1(\s)e_1+\cdots+ p_k(\s)e_k: \quad deg_{\s}(p_1(\s))\leq i \cdots
deg_{\s}(p_k(\s))
\leq i\}$$ 
being ${e_1,\cdots, e_k}$ a $B\{\s\}$-basis of ${\cal N}$.

\bigskip

\begin{defn}
An $A$-motive ${\cal N}$ of rank $n$ and $\s$-rank $k$ over $B$ is said pure if there exists
natural numbers
$  \bar{k}$, $\bar {n}$ such that  
$${\cal N}= {\underset {r\in {\Bbb N}}  \bigcup}{\cal N}_{r
\bar{n}}$$
 is a graded
$ {\underset {r\in {\Bbb N}} {\bigcup }H^0(X ,\o_X(r\bar{k} \infty) )  \otimes B} $-module. 
\end{defn}
In (2.3.1) we prove $n/k=\bar n /\bar k $. Recall that $\underset {r\in {\Bbb N}}
\bigcup H^0(\o_X(r \bar{k}\infty))=A$.

\bigskip

\begin{defn}\label{drin} A  $k$-elliptic sheaf, $({\cal E}_{jk},i_{jk},\tau)$, is a commutative
diagram of vector bundles of rank 
$n$ over $X\times S$, and injective morphisms of modules $\{i_h\}_{ h\in {\Bbb N}}$, $\tau$:

$$ \xymatrix { {\cdots} \ar[r] & 
{\cal E}_{-k} \ar[r]^{i_{-k}} & 
{\cal E}_{0} \ar[r]^{i_{0}} &
 {\cdots} \ar[r]^{i_{k(n-1)}}&
{\cal E}_{kn } \ar[r]^{i_{kn}} &
\cdots \\ 
\cdots \ar[ur]^\tau \ar[r] & 
F^* {\cal E}_{-k}\ar[ur]^\tau \ar[r]^{F^*i_{-k}} &
F^* {\cal E}_{0}\ar[ur]^\tau \ar[r]^{F^*i_{0}} &
{\cdots }\ar[ur]^\tau \ar[r] &
F^*{\cal E}_{ kn }\ar[ur]^\tau \ar[r]^{F^*i_{  kn }}&
{\cdots}}$$

satisfying:

a) For any $s\in S$, $deg(({\cal E}_{hk})_s)=n(g-1)+(h+1)k$.   $deg$ denotes the degree.  
 \bigskip \bigskip

b) For all $i\in \Bbb {\cal E}$, ${\cal E}_{(i+n)k}={\cal E}_{ik}(k\infty)(={\cal
E}_{ik}\underset {{\cal O}_{X\times S}}  
\otimes ({\cal O}_X(k\infty) \otimes {\cal O}_S))$. 
\bigskip \bigskip

c) ${\cal E}_{ik}+\tau(F^*{\cal E}_{ik})={\cal E}_{(i+1)k}$. 
\bigskip \bigskip

d) $R^1\pi_*({\cal E}_{-k})=0$. (As consequence  $\pi_*({\cal E}_0) $ is a locally
free $\o_S$-module  of rank $k$, in the case of elliptic sheaves this property is redundant).
  
 \end{defn}
 
  \bigskip

\bigskip

(2.3.1) Following \cite{D}, \cite{Mu}, we obtain an equivalence between the categories of pure
$A$-motives of rank $n$ and $\s$-rank $k$ and $k$-elliptic sheaves of rank $n$, just by
considering  ${\cal E}_{jk}=\widetilde{{\underset {r\in {\Bbb N}} 
\oplus}{\cal N}_{j+r\bar{n}}}$ as
$Proj({\underset {r\in {\Bbb N}}  \oplus}H^0(X,\o_X(r\bar{k}\infty)) \otimes B)$-graded module. $i_h$
(for each
$h\in {\Bbb Z}$) and $\tau$ are obtained from the natural inclusions among $\{{\cal
E}_{jk}\}_{j}$ and the multiplication on the left by $\s$ respectively.
  We can observe  that for $r>>0$  
$h^0({\cal E}_0((r+1)\bar{k}\infty))-h^0({\cal E}_0(r\bar{k}\infty))=rank_B({\cal
N}_{{(r+1)}\bar{n}}/{\cal N}_{{r }\bar{n}})=\bar{n}k$ moreover  using the Riemann-Roch theorem
$h^0({\cal E}_0((r+1)\bar{k}\infty))-h^0({\cal E}_0(r\bar{k}\infty))=n\bar{k}$ therefore
$n/k=\bar n/\bar k$. 
\bigskip

\begin{defn}
$E=\overset {k} { \underset {1 } \oplus }{{\Bbb G}_a} $ is said to be  an abelian $t$-module over
$A$ of rank
$n$ and $\s$-rank $k$ if $E$ is an $A$-module valued functor in the category of
${\Bbb F}_q$-schemes verifying:  $Hom_{End_S({\Bbb G}_a)_S}(E,({\Bbb G}_a)_S)$ is a
locally free
$A\otimes \cal O_S$-module of rank $n$, for each $S$. 
\end{defn}

\bigskip

Again,    paraphrasing \cite{An1} we know that the categories of abelian
$t$-modules over $A$ and $A$-motives are antiequivalent, where 
$${\cal N}=Hom_{End_S({\Bbb
G}_a)_S}(E,({\Bbb G}_a)_S)$$

 An $t$-module over $A$ is called pure when the associated $A$-motive is
it.

From now, we suppose $A$-motives and $t$-modules over $A$ are pure.
\bigskip

\begin{rem}\label{t}When $A={\Bbb F}_q[t]$ and $B$ is a perfect field $K$, 
 it can be added   another condition to the already  considered objects: there exists an element
$\theta \in K$ and $N \in {\Bbb N}$ such that 
$(t-\theta)^N.Lie(E)=0$, that is $(t-\theta)^N({\cal N}/\s {\cal N})=0$ for $A$-motives, so we
recover the usual definitions of $t$-motives and $t$-modules. Over an arbitrary ring $B$, this
condition means that the characteristic polinomyal of
$t.:{\cal N}/\s {\cal N} \to {\cal N}/\s {\cal N}$ is $(x^{p^j}-b)^m$ ($k=p^jm$ and
$(p,m)=1$) for some $b\in B$, $t.$ is the multiplication by $t$. 
\end{rem}

\bigskip
\bigskip

Now we are recalling the definitions of level structures for $A$-motives, $k$-elliptic sheaves and
$t$-modules over $A$. 

\begin{defn}
A twisted $\bar{D}$-level structure on a vector bundle of rank $n$ over $X\times S$ is a pair,  $(
M,f_{\bar{D}})$, with
$f_{\bar{D}}$  a surjective morphism of sheaves of ${\cal O}_{X\times S}$-modules         
$$f_{\bar{D}}:M\to ({{\Omega}_X}/{{\Omega}_X(-{\bar{D}})})^n  {\otimes}{{\cal O}}_S $$  
 \end{defn}

Two $\bar{D}$-level structures, $(M,f_{\bar{D}})$ and $(M',{f_{{\bar{D}}'}} )$, are said to be
equivalent  if there exists an isomorphism of sheaves of ${\cal O}_{X\times S}$-modules,
$\phi:M\to M'$, compatible with the twisted level structures. 
 
\begin{defn}
A twisted  ${\bar{D}}$-level structure, $({\cal E}_{jk},i_{jk},\tau, f_{\hat x})$, for a
$k$-elliptic sheaf of rank
$n$ over $ S$,  is a  twisted
${\bar{D}}$-level structure for each vector bundle $ {\cal E}_{jk}$ compatible with the morphisms
$\{i_{jk},\tau\}$.
 \end{defn}

\begin{defn}
A ${\bar{D}}$-level structure in a $A$-motive of rank $n$ and $\s$-rank $k$ over $B$, ${\cal N}$,
is a pair
$({\cal N}, h_{\bar{D}})$ where 

$$h_{{\bar{D}}}:{\cal N}\to    ({{\Omega}_A}/{{\Omega}_A(-{\bar{D}})})^n 
{\otimes}B$$  

is a surjective morphism of $A \otimes B\{\s\}$-modules.
$({{\Omega}_A}/{{\Omega}_A(-{\bar{D}})})^n   {\otimes}B$ is
$B\{\s\}$-module by defining $\s.(\omega \otimes b)=\omega \otimes b^q$
\end{defn}

\bigskip

\begin{defn}
If $E$ is a $t$-module over $A$ of rank $n$ and  $\s$-rank $k$ over $B$, a $\bar{D}$-module
structure for
$E$ is a pair $(E,i_{\bar{D}})$. $i_{\bar{D}}$ is an isomorphism of $A$-modules

$$i_{\bar{D}}:{E}_{I_{\bar{D}}}(B)\to \widetilde {(I_{\bar{D}}/A)^n}(B)$$

where ${E}_{I_{\bar{D}}}$ is the subgroup scheme of $I_{\bar{D}}$ division points of $E$ and 
$\widetilde {(I_{\bar{D}}/A)^n}$ is the constant sheaf of stalk $(I_{\bar{D}}/A)^n$. 
 \end{defn}
\bigskip

Remark  that    $A$-motives, $k$-elliptic sheaves and $t$-modules over $A$ have a level
structures on ${\bar{D}}$, if the restrictions to $I_{\bar{D}}\times S$ of $  {\cal N}/\s {\cal
N} $ and
${\cal E}_0/\tau(F^*{\cal E}_{-k})$ are
 zero and $\text{Ann}(lie(E))$ as $A\otimes \o_S$-module is such that
$\text{supp}(Ann(lie(E)))\cap (\text{supp}(\bar{D})\times S)=0$. 
\bigskip

\begin{rem}\label{niv}
One can check in a standard way that the moduli categories of $\bar{D}$-level structures for
$k$-elliptic sheaves and
$A$-motives  are equivalent. On other hand if $I_{\bar{D}}=(a)$ with $a\in A$, one obtains a
similar result for
$A$-motives and $t$-modules over $A$.
Following \cite{An1} we are going to recall how is settled this last "equivalence".

Let us consider  an $a$-level  structure (means $(a)_0$-level structure), $(E, i_a)$, on a
$t$-module over
$A$ , and
${\cal N}$ the  associated $A$-motive to the $t$-module over $A$, $E$, then $i_a$ gives an
$a$-level structure on
${\cal N}$:

let $\phi$ be the surjective morphism of $A\otimes B\{\s\}$-modules defined by
$$\phi:{\cal N}=B\{\s\}e_1\oplus \cdots \oplus B\{\s\}e_k \to Hom_{{\Bbb F}_q}({E}_{a}(B), B)$$

where $\phi(p(\s)e_i)(\alpha_1,\cdots,\alpha_k)=p(\s)(\alpha_i)$ and 
$$(\alpha_1,\cdots,\alpha_k)\in {E}_{a}(B)\subset E(B)={\Bbb G}_a(B)\oplus \cdots \oplus{\Bbb
G}_a(B)$$
then if $Spec(B)$ is connected  from the isomorphism of $A\otimes B$-modules \newline
$i_a:{E}_{a}(B)\overset {\sim} \to  ( a^{-1}A/A)^n$ we get 
$$ Hom_{{\Bbb F}_q}({E}_{a}(B), B)\overset {(i_a)_* } \simeq  Hom_{{\Bbb F}_q}((a^{-1}A/A)^n,
B)= Hom_{{\Bbb F}_q}((a^{-1}A/A)^n,{\Bbb F}_q )\otimes B$$
and   by using the residue pairing (for simplicity $  aA=m^r_x$ with $x\in X$)
$$ {\Omega}_A/ m^r_x{\Omega}_A\times {m^{-r}_x} /A \to k(x) \overset {tr} \to {\Bbb F}_q$$

 we obtain a $rx$-level structure for ${\cal N}$ in a direct way. $tr$ is the trace map in $x\in
Spec(A)$.

\end{rem}
\bigskip

Now we are going to sketch some facts about   Sato's infinite
Grassmannian and vector bundles. We can suppose $x$ rational to easy notation.

Let $\underline{{\Bbb F}_q[[t_x]]}$,  
$\underline{{\Bbb F}_q((t_x))}$ be the functors in the category of ${\Bbb F}_q$-schemes
associated to the ring and field of formal  series, namely 
$$\underline{{\Bbb
F}_q[[t_x]]}(S)= H^0(S)[[t_x]] \text{ and } 
\underline{{\Bbb F}_q((t_x))}(S)=  H^0(S)[[t_x]][{t_x}^{-1}]$$
 respectively, for each
scheme
$S$. We denote by $\o_S((t_x))$ and by $\o_S[[t_x]]$ to the sheaves in the category of $S$-schemes
$e^*(\underline{{\Bbb F}_q((t_x)) })$ and $e^*(\underline{{\Bbb F}_q[[t_x]] })$
 respectively, recall that $e$ is the structure map  $S \to {\Bbb F}_q $.

There exists a scheme 
$$Gr({\Bbb F}_q[[t_x]]^n\text{,}{\Bbb F}_q((t_x))^n)\text{ (Sato's infinite
Grassmannian scheme)}$$
 which points over an arbitrary scheme $S$, are    cuasicoherent ${\cal
O}_S$-submodules, 
$ {\cal L}$, of 
$\o_S((t_x))$, such that in the category of $S$-schemes is
 $$  \gamma^*{\cal L} \subset  \o_{S' }((t_x))^n$$
 where $\gamma:S'\to S$, 
also must exist some $r\in {\Bbb N}$ verifying
$$ {\cal L} +  {t_x}^r \o_{S }[[t_x]]^n  =\o_{S }((t_x))^n $$ 
and 
${\cal L} \cap {t_x}^r\o_{S }[[t_x]]^n $ is a coherent free $\o_S$-module (this
last rule if $S$ is compact). These submodules are called discrete submodules.

Let $\fu{\cal M}_{rx}^n$ be the moduli functor of vector bundles of rank $n$ over $X$ with level
structures on $rx$ and 
$$ {\fu{\cal M}}_{\hat{x}}^n={\underset {r\in {\Bbb N}}
 \varprojlim}\fu {\cal M}^n_{rx}$$
If one fix   an uniformizer for $x$, $t_x$,  we obtain by the relative Krichever's morphism,
 \cite{MR}, \cite{Q},  that
$ {\fu{\cal M}}_{\hat{x}}^n$ is the subfunctor of  the functor of points of 
 Sato's infinite grassmannian scheme, of discrete submodules which are
$H^0(X-\{x\},\o_X)$-submodules. Moreover, it is  representable by a closed subscheme of the
infinite grassmannian \cite{Al}.

The elements of $(M,f_{\hat{x}}) \in {\fu{\cal M}}_{\hat{x}}^n(S)$ are called formal
$x$-level structures, over
$X\times S$ in an analogous way one can define  twisted $x$-formal level structure. If $A_M$ is
the discrete submodule associated  to the formal
$x$-level structure over
$X\times S$, $(M,f_{\hat{x}})$, 
then $\pi_*( M(rx))=A_M \cap {t_x}^r\o_{S }[[t_x]]^n $ and 
$R^1\pi_*( M(rx))= \o_{S }((t_x))^n  /A_M + {t_x}^r\o_{S }[[t_x]]^n $

To see more details about these stamens, one can consult \cite{BeL}, \cite{Al}.

To study the moduli problem for $k$-elliptic sheaves we have to settle some definitions.

\begin{defn}\label{form}
A twisted $x$-formal level structure in a $k$-elliptic sheaf of rank $n$ over $S$,
$({\cal E}_{jk},i_{jk},\tau, f^{tw}_{\hat{x}})$, is a twisted $x$-formal level structure,
$f^{tw}_{\hat{x}} $in each 
$ {\cal E}_{jk}$, compatible with the morphisms $\{i_{jk},\tau\}$.
\end{defn}

Therefore, fixing an isomorphism of ${\Bbb F}_q[[t_x]]$-modules
$\mu_x:\Omega_x\overset {\sim}\to {\Bbb F}_q[[t_x]]$ we can obtain from a twisted $x$-formal level
structures a
$x$-formal level structure 
  for the  $k$-elliptic sheaves, $f_{\hat{x}} $, by
considering  $f_{\hat{x}} =(\overset n \oplus \mu_x)(f^{tw}_{\hat{x}})$.

\bigskip
By copying the results of \cite{Al},  we can obtain:

\begin{thm}\label{moduli}
The moduli functor, $\fu{ \cal D }_{\hat{x}}^{n,k}$, of $k$-elliptic sheaves of rank $n$ with
$x$-formal level structures, $({\cal E}_{jk},i_{jk},\tau, {f}_{\hat{x}})$, is representable by a
locally closed subscheme, $ { \cal D }_{\hat{x}}^{n,k}$, of
$Gr({\Bbb F}_q[[t_x]]^n\text{,}{\Bbb F}_q((t_x))^n)$.
\end{thm}

\begin{proof}
 A discrete submodule ${\cal L}_0$  gives a $k$-elliptic sheaf with formal level structure if it
verifies:

\begin{itemize}
 
\item  
${\cal L}_{0}$, is a $H^0(X-\{x\},\o_X)$-module.

\item  For all $1\leq i\leq n $, $(\sigma^*)^i{\cal L}_0\subseteq {\cal L}_{0}(k\infty)$ and
 ${\cal L}_0+\cdots+(\sigma^*)^i{\cal L}_0+\cdots+(\sigma^*)^n{\cal L}_0=
{\cal L}_{0}(k\infty)$.

\item For all $1\leq i\leq n $, ${\cal L}_{0}(k\infty)/{\cal
L}_0+\cdots+(\sigma^*)^{i }{\cal L}_0$ is locally free of rank $k(n-i)$.
Where ${\cal L}_{0}(k\infty)={\cal L}_{0}\underset {{\cal O}_{X-\{x\}}} \otimes {\cal
O}_{X-\{x\}}(k\infty)$

this gives a diagram of discrete subspaces:
$$ \xymatrix { {\cdots} \ar[r] & 
{\cal L}_{-k} \ar[r] & 
{\cal L}_{0} \ar[r] &
 {\cdots} \ar[r]&
{\cal L}_{nk } \ar[r] &
\cdots \\ 
\cdots \ar[ur] \ar[r] & 
\s^* {\cal L}_{-k}\ar[ur] \ar[r] &
\s^* {\cal L}_{0}\ar[ur] \ar[r] &
{\cdots }\ar[ur]  \ar[r] &
\s^*{\cal L}_{nk }\ar[ur]\ar[r]&
{\cdots}}$$
where the morphisms are the natural inclusions and $ {\cal L}_{ik} = {\cal L}_0 +\cdots
+(\s^*)^i{\cal L}_0$ if $i\geq 0$ and if $i\leq 0$ and $ik+\lambda kn\geq 0$, 
$ {\cal L}_{ik}={\cal L}_{ik+\lambda kn}(-\lambda k\infty)$. Recall that $\s^*{\cal L}_0$ is the
pull back
 of ${\cal L}_0$ by the  Frobenius morphism $\s$.

\item  If $({\cal E}_{-k},f_{\hat{x}})$ is the formal level associated to ${\cal L}_{-k} $,
$deg({\cal E}_{-k})=n(g-1)$ and $R^1\pi_*({\cal E}_{-k})=0$  therefore $\pi_*({\cal E}_0)$ is
locally  free of rank $k$. In the case of elliptic sheaves it is enough to impose the condition
over the $deg$.
 \end{itemize}

Conversely if we have $({\cal E}_{jk},i_{jk},\tau,f_{\hat{x}})$ over $S$, ${\cal L}_{ik}$ is 
$${\cal L}_{ik}=f_{\hat{x}}({{\cal E}_{ik}}_{\vert X-\{x\} \times S})\subset 
\o_S[[t_x]]^n_{\vert X-\{x\} \times S}=\o_S((t_x))^n $$
\end{proof}

When $A={\Bbb F}_q[t]$ the condition of the Remark \ref{t} is verified over a closed subscheme of
$ { {
\cal D }}_{\hat{x}}^{n,k}$.

Actually the $S$-points of $ { \cal D }_{\hat{x}}^{n,k}$ are not $A$-motives of $\s$-rank $k$ but
$A$-motives of locally (over $S$) $\s$-rank $k$, it depends on $\pi_*({\cal E}_0)$, in the case
of elliptic sheaves there are no problems with this because $\pi_*({\cal E}_0)$ is free of rank
$1$.

%%%%%%%%%%%%%%%%%%%%%%%%%%%%%%%%%%

\section{Subspace of uniformizers and locally dense subspaces for $k$-elliptic sheaves.}

We shall denote by
$ {\cal L}_0^{\cal E} $ to the discrete subspace ${\cal L}_0$ associated to  a $k$-elliptic sheaf,
 $({\cal E}_{jk},i_{jk},\tau,f_{\hat{x}})$,  with a $x$-formal
level structure.

\begin{defn}\label{den}
A cuasicoherent $\o_S$-submodule $D$ of
$\o_S[[t_x]]^n$, such that for each $S$-scheme $\gamma:S'\to S$,  $\gamma^*D\subset 
\o_{S'}[[t_x]]^n$, is said a locally dense subspace over $S$
if for some $m\geq 0$ and for every $h\geq
0$, $ D\cap {t_x}^{m+h}\o_S[[t_x]]^n/D\cap {t_x}^{m }\o_S[[t_x]]^n $ is  a  locally free
$\o_S$-module of rank $nh$.   
\end{defn}

\bigskip

\begin{propo}\label{subde}

If $({\cal E}_{jk},i_{jk},\tau,f_{\hat{x}})$ is a $k$-elliptic sheaf of rank $n$ with a $x$-formal
level structure over $S$, then its associated $A$-motive of $\s$-rank $k$ and rank $n$ with a
formal level structure on $x$, $({\cal N}, g_{\hat x})$ can be characterizated as a locally
dense  subspace,
$D_{\cal N}$, of 
$\o_S[[t_x]]^n $.
\end{propo}

\begin{proof} 
 The locally dense subspace, $D_{\cal N}$, is $g_{\hat x}({\cal N})$, where $g_{\hat x}$ is the
level structure
$g_{\hat x}:{\cal N}\to \o_S[[t_x]]^n $, ( we have fixed
$\mu_x:\Omega_x\overset {\sim}\to {\Bbb F}_q[[t_x]]$). Moreover, since ${\cal N}=\underset {j\geq
0}
\bigcup \pi_*({\cal E}_j)$,
$D_{\cal N}$ is $\underset {i\geq 0} \bigcup  ( {\cal L}_j^{\cal
E}\cap \o_S[[t_x]]^n ) $ and as $f_{\hat{x}}(\pi_*({\cal E}_j))={\cal L}_j^{\cal
E}\cap \o_S[[t_x]]^n $ we conclude. 
\end{proof}

\bigskip

Now we are going to characterize $k$-elliptic sheaves of rank $n$ with formal level structures as
$k$-dimensional vector subspaces of
$ {\Bbb F}_q[[t_x]]^n $. For this we need to precise the relative definition of these
subspaces:

\begin{defn}\label{grass}
$H$ is said a $k$-dimensional vector
subspaces of
$\o_S[[t_x]]^n $ over $S$ if $H$ is a  coherent locally free submodule of rank $k$ of
$\o_S[[t_x]]^n $, in the sense explained before. i.e: $\gamma^*H \subset \o_{S'}[[t_x]]^n$, for
each $S$-scheme   $S'$ ($\gamma:S'\to S$). 
\end{defn}
\bigskip

It is not hard to prove that there exists a grassmannian scheme (not finite dimensional),
$Gr^k({\Bbb F}_q[[t_x]]^n)$, for these objects. 

In the next theorem, we are going to show ${\cal D}_{\hat x}^{n,k}$ is a subscheme
of this Grassmannian.
\bigskip

\begin{thm}\label{D}

There exist an injective morphism of schemes

$$\psi:{\cal D}_{\hat x}^{n,k} \to Gr^k({\Bbb F}_q[[t_x]]^n).$$
\end{thm}

\begin{proof}
\rm{$\psi$ is defined in the functors of points of the schemes ${\cal D}_{\hat x}^{n,k}$  and
$Gr^k({\Bbb F}_q[[t_x]]^n)$
by 
$$\psi({\cal E}_{jk},i_{jk},\tau,f_{\hat{x}})={\cal L}_0^{\cal E}\cap\o_S[[t_x]]^n
=f_{\hat{x}}(\pi_*{\cal E}_0)=H$$
 this is a subspace in the sense of the
last definition, because of standard properties for discrete subspaces and vector bundles. $\psi$
is injective since $H$ generates the locally dense
subspace $D_{\cal N}$, as $A\otimes {\cal O}_S\{\s\}$-module.}\end{proof}

When $k=1$ and $n=1$ (i.e classical elliptic sheaves of rank $1$) $\text{Im}(\psi)$ takes valors on
the subfunctor in the category of ${\Bbb F}_q$-schemes
$H^0(S)[[t_x]]^*/H^0(S)^*\subset  Gr^1({\Bbb F}_q[[t_x]])(S)$, where $()^*$ means  units of
rings, because the characteristic of the elliptic sheaves is away from $x$.

 \begin {defn}\label{unif} 
  $s_1,\cdots,s_k \in H^0(S)[[t_x]]^n$ are uniformizers for a $k$-elliptic sheaf of rank
 $n$ (or $A$-motives or $t$-modules over $A$) over $S$ if it is a basis for the associated
$\o_S$-modulo
$H$.
 \end{defn}
For elliptic sheaves with   $x$-formal level structures  there  always exists uniformizers since 
 that $H$ is free:
$$H={\cal L}_0^{\cal E}\cap\o_S[[t_x]]^n/{\cal
L}_0^{\cal E}\cap \Gamma\o_S[[t_x]]^n={\o_S}$$
 because 
$${\cal L}_0^{\cal E}\cap \Gamma\o_S[[t_x]]^n=\o_S((t_x))^n/{\cal L}_0^{\cal E}+
\Gamma\o_S[[t_x]]^n=\o_S((t_x))^n/{\cal L}_0^{\cal E}+
 \o_S[[t_x]]^n=0$$
where $\Gamma$ is the diagonal matrix $(t_x,1,\cdots,1)$.
For general $k$-elliptic sheaves there just exists uniformizers locally over $S$.

\bigskip

\begin{rem}
For $A={\Bbb F}_q[t]$, $x=(t)$ and $n=k=1$ (Carlitz's modules), this last theorem has a connection
with Serre's class field theory:

If we denote  by  $J_{{\Bbb P}_1 \text{,}\hat x}$ (local Jacobian) the representant scheme of
the  functor 
$$S\to H^0(S)[[t_x]]^*/H^0(S)^*$$
   $J_{{\Bbb P}_1 \text{,}\hat x}$ is obviously 
an open  subscheme of 
$Gr^1({\Bbb F}_q[[t_x]] )$. Moreover, by Serre's class field theory
 $$\psi:{\cal D}_{\hat x}^{1,1} \to J_{{\Bbb P}_1 \text{,}\hat x}(\ref{D})$$  
is the blowing up of the
Albanese morphism 
$$Spec({\Bbb F}_q[t^{-1}])\to J_{{\Bbb P}_1 \text{,}\hat x}$$
 by the isogeny 
$$Id/F^*:J_{{\Bbb P}_1 \text{,}\hat x}\to J_{{\Bbb P}_1 \text{,}\hat x}$$
 here the Albanese
morphism is given by the formal  serie $(1-t^{-1}t_x)^{-1}$. Therefore, one universal
uniformizer $s_1$   (unique up to units) can be obtained of the relation $s_1/\s^*
s_1=(1-t^{-1}t_x)^{-1}$. The coefficients of $s_1$ are the roots of the  Carlitz's polinomyals
\cite{H}.  For  details about the local Jacobian, $J_{{\Bbb P}_1 \text{,}\hat x}$ to see
\cite{C},
\cite{KSU}, \cite{AMP}.

\end{rem}

\bigskip

\begin{rem}\label{com}
Following the Remark \ref{niv} we can do an explicit calculation of $H$ (here $x$ is not
necessarily rational), we suppose
$S=Spec(B)$ and connected. If $(E,i_{\hat x})$ is a $t$-module over $A$ of rank $n$ and $\s$-rank
$k$ with a
$x$-formal level structure: 
$$i_{\hat x}:E_{\hat x}(B)=
\underset {r\in {\Bbb N}} \varinjlim E_{m_x^r}(B)\overset {\sim} \to (k(x)((t_x))/k(x)[[t_x]])^n$$
and $({\cal N}=\underset {i=1}{
\overset n
\oplus}B\{\s\}e_i,g_{\hat x})$ the associated $A$-motive with a
$x$-formal level structure then  
$H=<s_1,\cdots,s_k>$, where
$s_i$ is the image of $e_i$ by the level structure morphism $g_{\hat x}$. By direct
calculation we get  
$$(s_i)_\xi=\underset {h\geq 0} \sum(\underset {1\leq c\leq deg(x)} \sum
\phi(e_i)(\alpha^\xi_{c,h})\omega_c t_{x}^{h})$$

(similar to \cite{An1} in Theorem 5)
where $\{\omega_c\}_{1\leq c\leq deg(x)}$ is a ${\Bbb F}_q$-orthonormal basis of $k(x)$
for the pairing trace. $(s_i)_\xi$ is the $\xi^{\text{th}}$-component of $s_i\in ((k(x)\otimes
B)[[t_x]])^n$ and $\alpha^\xi_{c,h}\in \underset {r\in {\Bbb N}} \varinjlim E_{m_x^r}(B)$
verifying
$$i_{\hat x}(\alpha^\xi_{c,h}) =(0,\cdots,\overset{\overset \xi { \smile}}
{\omega_ct_x^{-1-h}},\cdots,0)\in (k(x)((t_x))/k(x)[[t_x]])^n$$ 
 in this last equality we have used the fixed isomorphism of $k(x)[[t_x]]$-modules:
$\mu_x:\Omega_x\overset {\sim}\to k(x)[[t_x]]$.
\end{rem}
\bigskip \bigskip 

From now we are studying the case of classical elliptic sheaves in a more precise way. 
%%%%%%%%%%%%%%%%%%%%%%%%%%%%%%%%%%%%%%%%%%%%%%%%%%%%%%%%%%%%
\section{Discrete subspaces and elliptic Baker function.}

In this section we answer the following questions: How acts the classical arithmetic groups on
$H=<s_1>$? ( unfortunately, the classical arithmetic groups do not act on general $A$-motives,
because d) of the Definition 2.3). Is it  possible to describe basis for $D_{{\cal N}}$ and $ {\cal
L}_0^{\cal E} $ in terms of
$<s_1>$? Which is the counterpart of the Baker function in this context?

Let $\beta$ be an element of $Gl_n(k(x)((t_x)))$, with $v_x(det(\beta))=l$, then the discrete
subspace associated to the classical action, $ ({\cal E}_{j},i_{j},\tau,f_{\hat{x}})^\beta$, of
$\beta$, in the elliptic sheaves of rank $n$  with a
$x$-formal level structure over $S$
is: 
$\beta{\cal L}_{-l}^{\cal E}$ where the operation of $\beta$ in the discrete subspaces is
defined in the obvious way, so $<s_1^\beta>=\beta{\cal L}_{-l}^{\cal E}\cap
e^*(\underline{k(x)[[t_x]]^n})$. Then if $l\geq 0$, $\beta^{-1}s_1^\beta\in  {\cal
L}_{0}^{\cal E}$ because

$${\cal L}_{-l}^{\cal E}\cap e^*(\underline{\beta^{-1} (k(x)[[t_x]]^n)})\subseteq {\cal
L}_{0}^{\cal E}\cap e^*(\underline{\beta^{-1} (k(x)[[t_x]]^n)})$$ 

and if $l\leq 0$, $\beta^{-1}s_1^\beta\in D_{{\cal N}}$.

\bigskip
In the next theorem, we are going to get a basis for $ {\cal L}_{0}^{\cal E}$, by using 
the action of $Gl_n(k(x)((t_x)))$ on $<s_1>$. Again to do this we suppose $x$ rational.

\begin{lem}\label{base}
  $\{{\beta^{-1}_{(i,r)}}{s_1^{\beta_{(i,r)}}}\}_{r\geq 0,1\leq i\leq n}$ is a basis of $
{\cal L}_{0}^{\cal E}$ as $\o_S$-module, where
$\beta_{(i,r)}$ is the diagonal matrix $diag(1,\cdots,1,\overset   {\overset i
\smile} {t_x^{r}},1,\cdots,1)$.(The elliptic sheaf associated to $s_1$ is defined over $S$).
\end{lem}

\begin{proof}
By the before observation
 $$\beta^{-1}_{(i,r)}{s_1}^{\beta_{(i,r)}}\in {\cal L}_{0}^{\cal E}\cap
(\o_S[[t_x]] \oplus\cdots \oplus t_x^{-r}
\o_S[[t_x]]\oplus\cdots\oplus \o_S[[t_x]])$$
 moreover we deduce that 
 
$$\{\beta^{-1}_{(i,r)}{s_1}^{\beta_{(i,r)}}\}_{r\geq
0,1\leq i\leq n} $$ 
are linearly independent, since that  the $k^{th}$-components of
${s_1}^{\beta_{(i,r)}}$ are units in $\o_S[[t_x]]  $ because if $\beta\in
Gl_n({\Bbb F}_q((t_x)))$ with $deg(\beta)>0$ is $ {\cal L}_{0}^{{\cal E}^{\beta_{(i,r)}}}\cap
\beta \o_S[[t_x]]^n=0$. On the other hand, as 
$$  \underset {h\in {\Bbb N}} {\bigcup } ({\cal L}_{0}^{\cal E}\cap
t_x^{-h}\o_S[[t_x]]^n))={\cal L}_{0}^{\cal E}  $$ 
 and $\text{rank}_{\o_S}({\cal L}_{0}^{\cal E}\cap
t_x^{-h}\o_S[[t_x]]   ^n))=nh+1$ we conclude that 
$\{\beta^{-1}_{(i,r)}{s_1}^{\beta_{(i,r)}}\}_{r\geq 0,1\leq i\leq n} $ are generators, therefore it
is a basis. \end{proof}
\bigskip

\begin{lem}\label{dete} For ${r< 0,1\leq \xi\leq n}$ is 
$${s_1}^{\beta_{(\xi,r)}}=  \beta_{(\xi,r)} \text{det}\pmatrix  \alpha_0^\xi & \alpha_1^\xi &
\cdots &
\alpha_{-r-1}^\xi & s_1\\   
(\alpha_0^\xi)^q & (\alpha_1^\xi)^q & \cdots & (\alpha_{-r-1}^\xi)^q &
\sigma^*(s_1)\\ 
\cdot & \cdot & \cdots & \cdot & \cdot \\ 
\cdot & \cdot & \cdots & \cdot & \cdot \\
(\alpha_0^\xi)^{q^{-r  }} & (\alpha_1^\xi)^{q^{-r }} & \cdots &  
 (\alpha_{-r-1}^\xi)^{q^{-r } } & (\sigma ^*)^{-r }(s_1) 
    \endpmatrix $$
up to units of $H^0(S)$. We have used the notation of the Remark \ref{com}, recall that in this
case as
$k=1$, $\phi(e_1)=Id$. 
\end{lem}
\begin{proof} The proof follows 
$  <{s_1}^{\beta_{(\xi,r)}}>=\beta_{(\xi,r)}{\cal L}_{-r}^{\cal E}\cap
 \o_S[[t_x]]^n   $ and 
$$\beta_{(\xi,r)}\text{det}\pmatrix  \alpha_0^\xi & \alpha_1^\xi &
\cdots &
\alpha_{-r-1}^\xi & s_1\\   
(\alpha_0^\xi)^q & (\alpha_1^\xi)^q & \cdots & (\alpha_{-r-1}^\xi)^q &
\sigma^*(s_1)\\ 
\cdot & \cdot & \cdots & \cdot & \cdot \\ 
\cdot & \cdot & \cdots & \cdot & \cdot \\
(\alpha_0^\xi)^{q^{-r  }} & (\alpha_1^\xi)^{q^{-r }} & \cdots &  
 (\alpha_{-r-1}^\xi)^{q^{-r } } & (\sigma ^*)^{-r }(s_1) 
    \endpmatrix \in \beta_{(\xi,r)}{\cal L}_{-r}^{\cal E}   $$
\end{proof}
Note that if $d \in {\Bbb N}$ is such that $m_\infty^d=bA_x$ is principal and 
$$H^0(X-\{x\}, \o_X)=A_x\subset {\Bbb F}_q((t_x))$$
then $b=ut^{-d}_x$ where
$u$ is a unit in  ${\Bbb F}_q((t_x))$, therefore ${s_1}^{diag(t^{-d}_x,\cdots,t^{-d}_x)}=u^{-1}s_1$
since that
$s_1^{ut^{-d}_x}=s_1$ (up to units).

This remark  together Lemma \ref{base} and
$$\beta_{(\xi,1)}=t^d_x\beta_{(1,-d)}\cdots \beta_{(\xi,-d+1)} \cdots \beta_{(n,-d)} $$
 gives us a
explicit basis of ${\cal L}_0^{\cal E}$  ($s_1^{ t^{d}_x}=us_1$), therefore  $H=<s_1>$ determines
${\cal L}_0^{\cal E}$.
\bigskip\bigskip

The Proposition 5.1 of \cite{SW} about the Baker function gives here an analogous result. Before
we are going to state some notation

$ \underline {Gl_n^0} ({\Bbb F}_q((t_x)))(S)=\{g\in Gl_n (H^0(S)^\sigma((t_x)))/ \text{ with
}v_x(det(g)_s)=0 \text{ for all } s\in Spec(H^0(S))\}$  where $H^0(S)^\sigma=\{c\in H^0(S)/
c^q=c\}$.

\begin{propo} For each ${\cal L}_0^{\cal E}$, with the usual notation in this paper, there exists
a unique, up to units of
$H^0(S)$, (algebraic) function:
$$ \Psi_{{\cal L}_0^{\cal E}}:\underline {Gl_n^0}({\Bbb F}_q((t_x)))(S)\to H^0(S)((t_x))^n$$
such that:

1) $\Psi_{{\cal L}_0^{\cal E}}(g)\in {\cal L}_0^{\cal E}$ for every $ g \in \underline
{Gl_n^0}({\Bbb F}_q((t_x)))(S)$.

2) $\Psi_{{\cal L}_0^{\cal E}}(g)=g^{-1}( \underset {i=0}{\overset {\infty} \sum} a_i(g)t_x^i)$
with $a_i$ (algebraic) functions from $\underline {Gl_n^0}({\Bbb F}_q((t_x)))(S)$ to $H^0(S)^n$,
$a_0$ takes valors on ${H^0(S)^*}^n$.
\end{propo}

\begin{proof}
This function is defined by $ \Psi_{{\cal L}_0^{\cal E}}(g)=g^{-1}s_1^g$. It is unique by the
properties of the elliptic sheaves. Obviously it extends on the whole of $\underline
{Gl_n }({\Bbb F}_q((t_x)))(S)$.
\end{proof}
This function is algebraic since it can be interpreted as a morphism of functors in the category
of $S$-schemes.

From the Lemma \ref{dete} we can get an explicit calculation for $\Psi_{{\cal L}_0^{\cal E}}(g)$,
with
$g\in {Gl_n }({\Bbb F}_q((t_x)))$, since diagonal matrixs and ${Gl_n }({\Bbb F}_q[[t_x]])$
generate
${Gl_n }({\Bbb F}_q((t_x)))$ and $ \Psi_{{\cal L}_0^{\cal E}}(g)=s_1$ for all $g\in {Gl_n }({\Bbb
F}_q[[t_x]])$.

In an obvious way we can translate these results to Krichever $D$-modules \cite{D}.
%%%%%%%%%%%%%%%%%%%%%%%%%%%%%%%%%%%%%%%%%%%%%%%%%%%%%%%%%%%%%%%%%%%%%%%%%%%%%%%%

%\bigskip

%\begin{rem}
%With the behaviour of the action of  $Gl_n({\Bbb F}_q[[t_x]]) $ in
%${\cal D}^{n,1}_{\hat x}$, we can get a connection with the Fourier-Mukai transform : If $(L,
%h_{\hat x})$ is a line bundle of degree $0$ on $X$, with a $x$-formal level structure and
%$(E_i,i_j,\tau,f_{\hat x})$ is a universal elliptic sheaf with a $x$-formal level structure and
%its subspace of uniformizers $<{S}_1>$, then $<{S}_1>^{(L, h_{\hat
%x})}=(\pi_2)_*((\pi_1)^*(L) {\otimes}_{\o_{X\times {\cal D}^{n,1}_{\hat x}}} E_0) $ where
%$<{S}_1>^{(L, h_{\hat x})}$ is the action of $(L, h_{\hat x})$ on the $
%{\o}_{{\cal D}^{n,1}_{\hat x}}$-subspace $<{S}_1>$, explained before. Here the $\pi_i$ are the
%natural projections
%$$\xymatrix { & {X\times {\cal D}^{n,1}_{\hat x}}\ar[dl]_{\pi_1} \ar[dr]^{\pi_2}& \\
%X & & {{\cal D}^{n,1}_{\hat x}}}$$ 
%\end{rem}

%\bigskip \bigskip
%%%%%%%%%%%%%%%%%%%%%%%%%%%%%%%%%%%%%%%%%%%%%%%%%%%%%%%%%%%%%%%%%%%%%%%%%%%%%%%%

\section{Elliptic sheaf triples.}

In this section, we are going to translate
the results of the  antiequivalence of the super Krichever functor, \cite{MR}, \cite{SW}, into our
setting of 1-dimensional subspaces of uniformizers:

Paraphrasing  the classic case for Schur pairs, data $(\infty, X, x, g_{\hat x}, ({\cal
E}_{},i_{j},\tau,f_{\hat{x}}))$ have associated triples $(\infty, A_x, {\cal L}_0^{\cal E})$,
where
$x$ is a rational point of the nonsingular projective curve $X$, $\infty$ is a point of
$X-\{x\}$, $({\cal E}_{j},i_{j},\tau,f_{\hat{x}})$ an elliptic sheaf of rank $n$ on
$A=H^0(X-\{\infty\},\o_X)$ with a
$x$-formal level structure.
$g_{\hat x}:\o_X\to  {\Bbb F}_q[[t_x]]$ a $x$-formal level structure on $\o_X$ which is a ring
morphism, compatible with $f_{\hat {x}}$ and the $\o_x$-module structure of ${\cal E}_{j}$,  $A_x$
is the subspace discrete of
$ {\Bbb F}_q((t_x))$ defined by
$A_x=g_{\hat x}(H^0(X-\{x\},\o_X))$. ${\cal L}_0^{\cal E}$ is the discrete subspace associated
to 
$({\cal E}_{j},i_{j},\tau,f_{\hat{x}})$ in the last section. We name elliptic sheaf triples to  
triples $(\infty, A_x, {\cal L}_0^{\cal E})$.
 
Bearing in mind   \cite{MR}, \cite{SW},
 it is not hard to set up an one to one correspondence between isomorphisms of data 
 $$(\infty, X,x,g_{\hat x}, ({\cal E}_{j},i_{j},\tau,f_{\hat{x}}))$$
and isomorphisms of  elliptic sheaf triples. 

\begin{rem}\label{Schur} If we fix $d\in {\Bbb N}$ with $m_\infty^d=bA_x$, $<s_1>$ determines
$(\infty, A_x, {\cal L}_0^{\cal E})$: because
$<s_1>$ and the action of $Gl_n({\Bbb F}_q((t_x)))$ gives   a  basis of ${\cal
L}_0^{\cal E}$ (Lemma \ref{dete}) moreover since
$A_x$ is non singular 
$$A_x=\{c\in 
 {\Bbb F}_q((t_x))/ c.{\cal L}_0^{\cal E}\subseteq{\cal L}_0^{\cal E}\}$$ 
and $\infty=supp({\cal L}_0^{\cal E}+\s^*({\cal L}_0^{\cal E})/{\cal L}_0^{\cal E})$.  
$d$ is said to be the exponent of the elliptic sheaf triple $(\infty, A_x, {\cal L}_0^{\cal E})$.
Theorem \ref{D} shows this result by assuming  $\infty$ and $A_x$ fixed.
\end{rem}

We are going to use this last result  to study the relative  elliptic
sheaf triples  and its moduli problem.

Let $f:S\to \bar{S}$ be a morphism of shemes with $\bar{S}$ a scheme where the Frobenius morphism
$\sigma$ is the identity.

\begin{defn}\label{tri}
An elliptic sheaf triple of rank $n$ and exponent $d$ over $f:S\to \bar{S}$, is a triple 
$({ I}, { A}_x, {\cal L}_0)$ where ${ A}_x$ is a ring discrete subspace  of
$\o_{\bar{S}}((t_x)) $ together a discrete ideal
${ I}$ of $\o_{\bar{S}}((t_x)) $ such that ${ A}_x/{ I}=\o_{\bar{S}}$, $I/I^2$ is locally free of
rank $1$ over
$\o_{\bar{S}}$ and ${  f^*I}^d=ut_x^{-d}f^*{A}_x$ with
$u$ a unit in $H^0(S)^\sigma[[t_x]] $  $(H^0(S)^\sigma=\{c\in H^0(S)/ c^q=c\})$.
${\cal L}_0$ is a discrete subspace of $\o_S((t_x))^n$ which is an $f^*{
A}_x$-module fulfilling the following conditions :

1) $f^*I( {\cal L}_0+ \s^*{\cal L}_0+\cdots+ (\s^*)^n{\cal L}_0)\subseteq  {\cal
L}_0$.

2) For every $1\leq i \leq n$ 
$$ {\cal L}_0/ f^*I({\cal
L}_0+ \s^*{\cal L}_0+\cdots+ (\s^*)^{i}{\cal L}_0) $$
 is locally free of rank $ n-i $.

3) For every $s\in S$ 
$$dim_{k(s)}({{\cal L}_0}_s\cap
k(s)[[t_x]]^n)-dim_{k(s)}k(s)((t_x)))^n/{{\cal L}_0}_s+ k(s)[[t_x]]^n =1.$$

4) $(A_x)_s$ is a Dedekind ring for each $s\in \bar{S}$.

We denote by $F^d_n$ to the functor of elliptic sheaf triples of rank $n$ and exponent $d$.
\end{defn}

\bigskip
As in the last section, $({ I}, { A}_x, {\cal L}_0)$  has associated a  $1$-dimensional
subspace of uniformizers, $H=<s_1>={{\cal L}_0} \cap
\o_S[[t_x]]^n     $, futhermore this subspace  determines the triple $({ I}, { A}_x, {\cal L}_0)$:

\begin{propo}
If $({ I}, { A}_x, {\cal L}_0)$ and $({ I}', {{ A}'}_x, {{\cal L}'}_0)$ are two
elliptic sheaf triples of rank $n$ and exponent $d$, over a scheme $S$ with subspace of
uniformizers
$H=<s_1>$,
 then
$({ I}, { A}_x, {\cal L}_0)=({ I}', {{ A}'}_x, {{\cal L}'}_0)$.
\end{propo}
\begin{proof}
Since $k(s)={\Bbb F}_q$ for each $s\in \bar{S}$, $(\sigma=Id)$, Remark \ref{Schur} shows that 
$$({ I},
{ A}_x, {\cal L}_0)_{\vert f^{-1}(s)}=({ I}', {{ A'}_x}, {{\cal L}'}_0)_{\vert f^{-1}(s)}$$
 where $f$
is the
 morphism from $S$ to $\bar{S}$, thus the closed subscheme $C$ where
${ {\cal L}'}_0= {\cal L}_0$ (\cite{Al} Theorem 3.7) 
contains $f^{-1}(s)$ for each $s\in \bar{S}$, we have to prove that $C=S$. 
Locally $S=Spec(B)$,  $C=Spec(B/J)$  and  $J\subseteq \underset {s} \bigcap
m_sB$, with $m_s$
  the prime ideal associated to $s\in \bar{S}$,  bearing in mind that $\sigma$ is the identity on
$\bar{S}$ is
$J_{m_s}=0$ for each
$s\in \bar{S}$, so $J=0$.

On the other hand a similar argument proves that, ${ I}={ I}'$ and ${ A}_x= {{
A'}_x}$. 
\end{proof}
This last proposition gives an injective morphism of functors between $F^d_n$ and the functor of
points of
$Gr^1(({\Bbb F}_q[[t_x]])^n)$.
\bigskip

\begin{thm} $F^d_n$ is representable.
\end{thm}
\begin{proof}
We are going to sketch the proof.  Bearing in mind Sato's infinite grassmannian scheme and
(\cite{Al} Theorem 3.7), the set of points which define discrete
subspaces $(L' ,L)$ verifying $L.L=L$, ${\Bbb F}_q\subset L$ and $L'.L=L'$ is a subscheme of 
$$Gr({\Bbb
F}_q[[t_x]]\text{,}{\Bbb F}_q((t_x)))\times Gr({\Bbb F}_q[[t_x]]\text{,}{\Bbb F}_q((t_x)))$$
 so the functor of pairs $(I,A_x)$ 
 where
$A_x$ is a discrete subring over $S$ and $I$ is an discrete ideal of $A_x$ with the added
conditions
$A_x/I=\o_S$ and $({A_x})_s$ is a Dedekind ring for   each $s\in S$ is representable by a scheme $Z$.

Let
$Z_\s$ be the subsheme of
$Z$ where the Frobenius morphism  $\s$ is the identity.

Let ${{\mathfrak L} }_0$ and
$({\cal I},{\cal A}_x)$ be the pull back of universal objects for discrete subspaces and
pairs $\{(I,A_x)\}$ over $ Z_\s$, respectively,   on the scheme
$$W=Gr({\Bbb F}_q[[t_x]]^n\text{,}{\Bbb F}_q((t_x))^n) \times Z_\s$$
 by the natural projections,  then to get a representant  for $F^d_n$ we must impose to
${{\mathfrak L} }_0$ and
$({\cal I},{\cal A}_x)$ similar
conditions to the Definition \ref{tri} ones: 

\begin{itemize}

\item ${{\mathfrak L} }_0$ is an ${\cal A}_x$-module.

\item  For every $s\in W$ 
$dim_{k(s)}({{\mathfrak L}_0}_s\cap (k(s)[[t_x]]^n)=1$ and 
$$dim_{ k(s)}k(s)((t_x)))^n/{{\mathfrak
L}_0}_s+ (k(s)[[t_x]])^n)=0.$$

\item If locally $<s_1>={{\mathfrak L}_0}\cap (H^0(W')[[t_x]])^n$, $s_1^{t_x^{-d}}=u^{-1}s_1$
where
$u$ is a unit in $H^0(W')^\s[[t_x]]$, being $W'$ the scheme where the last two conditions are
verified, the definition of $s_1^{t_x^{-d}}$ is obtained from the Lemma \ref{dete}.

\item ${\cal I}^d=ut_x^{-d}{\cal A}_x$.
 
\item ${\cal I}/{\cal I}^2$ is locally free of rank $1$.

\item   $ {\cal I}({\mathfrak L}_0+ \s^*{\mathfrak L}_0+\cdots+ (\s^*)^n {\mathfrak
L}_0)\subseteq  {\mathfrak L}_0$.

\item For every $1\leq i \leq n$ 
$$ {\mathfrak L}_0/{\cal I}({\mathfrak L}_0+
\s^*{\mathfrak L}_0+\cdots+ (\s^*)^{i} {\mathfrak L}_0)$$
 is locally free of rank $n-i$.

\end{itemize}

\end{proof}

%%%%%%%%%%%%%%%%%%%%%%%%%%%%%%%%%%%%%%%%%%%%%%%%%%%%%%%%%%%%%%%%%%%%%%%%%%%%%%%%

\section{Determinants for elliptic sheaves.}

In this section we are going to study the behaviour of the determinant for  discrete
and locally dense subspaces associated to elliptic sheaves. In our
setting, for discrete subspaces, we obtain the same results that \cite{An1}.

 One can eaily show that the determinant for an  elliptic sheaf of rank $n$,  with the obvious
definition of  determinant, is an elliptic sheaf of rank $1$,  but when we work with level
structures ,we have to twist the natural definition of
  determinant for elliptic sheaves with $x$-formal level structures  to get again elliptic
sheaves with $x$-formal level structures:

Recall from  definitions that if $({\cal E}_{j  },i_{j  },\tau,f_{\hat{x}})$ is an  elliptic
sheaf of rank
$n$ over
$S$ with a
$x$-formal level structure, then $f_{\hat{x}}$ is $(\overset n \oplus \mu_x) (
 f^{tw}_{\hat{x}})$ where 
$$f^{tw}_{\hat{x}}: {\cal E}_{0 } \to \underset
{m\in {\Bbb N}} {\varprojlim}(\Omega_X/\Omega_X(-mx))^n$$ 
is a twisted $x$-formal level structure on the elliptic sheaf and $ {\mu_x}$ is a fixed
isomorphism from $ \Omega_x $ to ${\Bbb F}_q[[t_x]]$.

\begin{defn}\label{det}
We define the  elliptic sheaf determinant with
$x$-formal level structure, $det({\cal E}_{j },i_{j },\tau,f_{\hat{x}})$, of    an  elliptic
sheaf of rank $n$  over
$S$  with a $x$-formal level structure, $({\cal E}_{j },i_{j },\tau,f_{\hat{x}})$. As 
$$((\underset {\o_X\otimes \o_S} {\overset {n} {\bigwedge}} {\cal E}_{j }) \underset {\o_X} \otimes
\Omega_X^{\otimes1-n},\overset {n} {\wedge}i_j,\overset {n} {\wedge}\tau, \mu_x(
 f_{\hat{x}}'))$$
being $ f_{\hat{x}}'$ the formal level structure obtained from 
$$\overset {n} {\bigwedge}f^{tw}_{\hat{x}}: \underset {\o_X\otimes \o_S} {\overset {n} {\bigwedge}}
{\cal E}_{0 }  \to
\underset {m\in {\Bbb N}} {\varprojlim} \Omega^{\otimes n}_X/\Omega^{\otimes n}_X(-mx) $$

by tensoring by $\Omega_X^{\otimes 1-n}$. Recall that
a level structure over
${\cal E}_0$ (with some extra conditions) determines the level structure in the $k$-elliptic
sheaf, $({\cal E}_{j },i_{j },\tau,f_{\hat{x}})$, since characteristic is away from $x$.
\end{defn}

From this definition we can obtain the analogous result of \cite{An1}
for discrete subspaces  in a straightforward way:
\begin{propo}
If $ {\cal L}_0^{\cal E} $ is the associated discrete subspace to  $({\cal
E}_{j },i_{j },\tau,f_{\hat{x}})$ then 
  $$(\underset {\o_{X-\{x\}}\otimes \o_S} {\overset {n} {\bigwedge}}  {\cal L}_0^{\cal E})
\underset {\o_{X-\{x\}} }{\otimes} {\mu_x}( \Omega_{X-\{x\}})^{\otimes
1-n}((n-1)(g-1)\infty)$$
 is the discrete
subspace for  $det({\cal E}_{j },i_{j },\tau,f_{\hat{x}})$.
\end{propo}
\bigskip

Now we are   describing the determinant for locally dense subspaces, $D_{\cal N}$, associated to
elliptic sheaves. With the notation of the last definition we must to study the $A\otimes
\o_S\{\s\}$-module
$$\mu_x({f}_{\hat{x}}')(((\underset {\o_X\otimes \o_S} {\overset {n} {\bigwedge}} {\cal E}_0)\underset
{\o_X} \otimes
\Omega_X^{\otimes  1-n })_{Spec(A)\times S})=( \underset {A\otimes \o_S} {\overset {n}
{\bigwedge}} D_{\cal N}) \underset {A}
\otimes {\mu_x}(\Omega_A)^{\otimes  1-n }\subset \o_S[[t_x]]^n $$
  for   it we are going to consider $\underset {A\otimes \o_S} {\overset {n} {\bigwedge}}D_{\cal N}$.

\bigskip

\begin{lem}
If $\o_S\{\s\}(={\cal N}) $ is the $A$-motive associated to an elliptic sheaf of rank $n$  then
$\underset {A\otimes \o_S}{\overset {n} {\bigwedge}}\o_S\{\s\}$ is generated by $1\wedge
\s\wedge\cdots\wedge\s^{n-1}$ as $A\otimes\o_S\{\s\}$-module if and only if $\{1,
\s,\cdots,\s^{n-1}\}$ are generators of the $\o_S$-module $ \o_S\{\s\}/\o_S\{\s\}J$, where $J$ is
the ideal characteristic for the elliptic sheaf. ( $J=\text{Ann}({\cal E}_1/\tau (F^*{\cal E}_0))$
and $\o_X\otimes\o_S/J\simeq
\o_S$, with the standard notation.
\end{lem}

\begin{proof}
\rm{To prove the lemma it is suffices to suppose that $S=Spec(K)$ is a field.  

$\underset {A\otimes K}{\overset {n} {\bigwedge}}K\{\s\}$ is
generated by $1\wedge \s\wedge\cdots\wedge\s^{n-1}$ as $A\otimes K\{\s\}$-module if and only the set of
points of $Spec(A\otimes K)$ where some of the subsets   
$$\{1,\s,\cdots,\s^{n-1}\},
\{\s,\s^2,\cdots,\s^{n}\} ,\cdots, \{\s^m,\s^{m+1},\cdots,\s^{m+n-1}\},\cdots$$
are a $A\otimes K $-basis
of $K\{\s\}$ is
$Spec(A\otimes K)$. If we denote by $D_m$ the effective divisor on $X_K$ where
$\{\s^m,\cdots,\s^{m+n-1}\}$ is not a basis, we have to show that $\underset {m} \bigcap
\text{supp}(D_m )=\phi$ if and only if $q$ is not in $\text{supp}(D_0)$, where $q\in
X_K$ is the characteristic of the elliptic sheaf. Since $D_{m+1}=F^*D_{m}+q$ and therefore
$D_{m+1}=F^*D_{0}+q+F^*q+\cdots+(F^*)^mq$ to conclude it is enough to prove $\underset {m} \bigcap
\text{supp}((F^*)^mD_0 )=\phi$.
 But if $y\in \underset {m} \bigcap
\text{supp}((F^*)^mD_0 )$ then there exists $z\in Spec(A)$ such that $\{1,\s,\cdots,\s^{n-1}\}$
are not linearly independent in $K\{\s\}/K\{\s\}m_z$ as $A/m_z\otimes K$-module. But
by induction over $r\leq n-1$ one can check that $ \{1,\s,\cdots,\s^{r}\}$ are linearly
independent for all $r\leq n-1$ in $K\{\s\}/K\{\s\}m_z$ as $A/m_z\otimes K$-module.}
\end{proof}

\bigskip
One can observe that in an elliptic sheaf  over $S$, the set of points where the lemma is
verified is an open subset of $S$. One can also check through  a similar reasoning as in
the foregoing lemma that  if the characteristic of an elliptic sheaf is a rational point of
$Spec(A)$, this elliptic sheaf verifies the conditions of the lemma. So we can hope that the
open subset  of  $ {{ \cal D}}_{\hat{x}}^{n}$ where the condition of the lemma is verified is
not empty.

 From these lemmas we obtain:
\begin{thm}
If $({\cal E}_j,i_j,\tau,f_{\hat{x}})$ is an elliptic sheaf of rank $n$ over $S$ in the
conditions of the previous lemma, with a
$x$-formal level structure and locally dense subspace $D_{\cal N}$ then  
$\overset {n} {\bigwedge}(D_{\cal N})$ is generated as $A\otimes \o_S\{\s\}$-module by the
determinant of the ''scattering'' matrix 
$$(s_1,\s^*s_1,\cdots,(\s^*)^{n-1}s_1)$$
 where $<s_1>$ is the subspace of uniformizers associated to $({\cal
E}_j,i_j,\tau,f_{\hat{x}})$   in the last section.
\end{thm}

\begin{rem}
It would be very interesting to obtain the subspace of  uniformizers associated to 
$det({\cal E}_j,i_j,\tau,f_{\hat{x}})$ in terms of $<s_1> \leftrightarrow
({\cal E}_j,i_j,\tau,f_{\hat{x}})$. When
$X={\Bbb P}_1$ this subspace   is $<det(s_1,\s^*s_1,\cdots,(\s^*)^{n-1}s_1)>$  because
$( 1,\s ,\cdots,\s^{n-1})$ is a ${\Bbb F}_q[t]\otimes \o_S$-basis of $\o_S\{\s\}$. Unfortunately, one
can check that
  this happens only when
$g=0$. 
\end{rem}

\bigskip \bigskip

{\bf{ Acknowledgment}} I would like to express my gratitude to my thesis advisor Proffesor J.M.
Mu\~noz Porras for his support and encouragement. I would also like to thank to Proffesor Pedro Arias
for useful conservations.

%%%%%%%%%%%%%%%%%%

\vskip1.5truecm { \'Alvarez V\'azquez, Arturo}\newline {\it
e-mail: } aalvarez@gugu.usal.es
\end{document}